\documentclass[10pt,english]{article}
\usepackage{mathptmx}
\usepackage[T1]{fontenc}
\usepackage[utf8]{inputenc}
\usepackage[provide=*,english]{babel}
\usepackage{amsmath}
\usepackage{amssymb}
\usepackage{amsthm}
\numberwithin{equation}{section}
\renewenvironment{cases}
 {\left\{\begin{array}{@{}l@{\qquad}l@{}}}
 {\end{array}\right.}
\usepackage{xurl}
\usepackage{setspace}
\usepackage{graphicx}
\usepackage[pass]{geometry}
\usepackage{afterpage}
\makeatletter
\newcommand{\restoretextgeometry}{%
  \restoregeometry
  \global\@colht=\textheight
  \global\@colroom=\textheight
  \global\vsize=\textheight
}
\makeatother
\AtBeginDocument{
  \setlength{\abovedisplayskip}{12pt}
  \setlength{\belowdisplayskip}{12pt}
  \setlength{\abovedisplayshortskip}{12pt}
  \setlength{\belowdisplayshortskip}{12pt}
}
\newtheorem{theorem}{Theorem}

\usepackage[pdfusetitle,
 bookmarks=true,bookmarksnumbered=false,bookmarksopen=false,
 breaklinks=true,pdfborder={0 0 1},backref=false,colorlinks=false]
 {hyperref}
\begin{document}
\title{Cumulant Generating Function of the Hill-Lane-Sudderth urn model}
\author{Simone Franchini{\normalsize\thanks{Correspondence: simone.franchini@istc.cnr.it}\thanks{CNR-ISTC, Via Gian Domenico Romagnosi 18, 00196 Rome, Italy}}}
\date{~}
\maketitle

\begin{abstract}
\setstretch{1.1}
We consider the Hill, Lane and Sudderth urn, in which a black or
white ball is added at each step with probabilities depending on
the current proportion of black balls. We solve the equation for
the limiting normalized cumulant generating function when the urn
function is invertible on the outer branches beyond the fixed
points and the corresponding inverses are Lipschitz. The solution
has a parametric representation through three quadratures of the
urn function and a limit of compositions of radicals, with no
derivatives or inverses of the urn function in the final formula.
We prove an error bound that does not accumulate along the
partition and remains valid when starting from the initial
singularity. Four examples with exact solutions illustrate the
two orientations of the construction, the non-affine case and
the presence of zeros.
\end{abstract}

\newpage

\section{Introduction}
\label{sec:introduction}

The Hill, Lane and Sudderth (HLS) urn is a two-color reinforcement
model: one ball is added at each step, and the probabilities of
the two colors depend on the current composition~\cite{HLS1980}.
The current proportion thus summarizes the past additions, while
the urn function determines how this history influences the next
addition. Varying this response gives a common mathematical
description of growth processes arising in different settings
\cite{Pemantle2007,Franchini_HLS_review}.

Among these processes are leaf counts in random trees~\cite{Bryc2009}
and elephant random walks with either one or a fixed number of draws
with replacement from the past~\cite{BaurBertoin2016,PodderRoy2026}.
Nonlinear responses arise naturally when past choices affect subsequent
preferences. In technological diffusion, current market shares influence
later adoptions through positive or negative feedback
\cite{DosiErmolievKaniovski1994}; increasing returns and social influence
can produce competing long-run outcomes~\cite{FranchiniBalzan2023}.
When return weights are proportional to a common power of wealth,
the two-agent model with wages~\cite{GottfriedGrosskinsky2025} also
has an autonomous urn description. Each allocation of the capital
return is an urn addition, with a probability determined by the
current wealth shares.
Deterministic wages make accumulated wealth an affine transformation
of the urn count, so its fluctuations can be studied through the same
process~\cite{Franchini_HLS_review}.

In irreversible magnetic growth, spins are added with probabilities
determined by the accumulated magnetization
\cite{KlymkoGarrahanWhitelam2017}. A nonlinear response can sustain
competing magnetizations, making the distribution of growth histories
relevant even when the fixed points are known. Related competition
occurs between two growing neuronal extensions, or neurites
\cite{KhaninKhanin2001}.
P\'olya urn models have also been used to describe distributions of
enantiomeric excess in asymmetric autocatalysis~\cite{Barabas2009}.
The review~\cite{Franchini_HLS_review} develops cooperative nonlinear
production rules as further illustrative urn models with several
possible final compositions.

Nonlinear urns also describe the accumulation of information.
In broadcasting on growing recursive graphs, a new vertex samples
a fixed odd number of earlier vertices, at least three, uniformly
with replacement. Each sampled color is independently flipped with
a probability strictly between zero and one half, and the new color
is chosen by majority vote
\cite{BriendDevroyeLugosi2025}. The probability of choosing a given color
is therefore a polynomial of its current fraction, and the urn records
the number of vertices of each color. In reinforced recommendation, the reduction
to one user and two items gives an urn with power-law reinforcement
\cite{Bellina2023,Franchini_HLS_review}. Related feedback rules govern
which new observations are added to a stored memory, as discussed in
the review. A different response arises in
urns with triggering: discovering an element expands the set of
available elements~\cite{Tria2014}. This mechanism leads to a rational
limiting urn response for the count of distinct discoveries, such as
new words in a text.
The comparison of history probabilities in~\cite{Franchini_HLS_review}
justifies transferring the limiting free energy to this autonomous
HLS description.

For random walks and polymers, the accumulated observable can also be
geometric. The range counts distinct visited sites; at a fixed length
it determines the number of steps returning to old sites.
Conditional continuation probabilities reproduce this return-count distribution,
linking urns to the Stanley polymer model~\cite{FranchiniBalzan2018}.
The same range statistics determine survival in the Rosenstock problem
with random traps, while the zero-return sector selects
self-avoiding walks~\cite{Franchini_HLS_review}. Attractive interactions
require the additional contact statistic, as in the study of interacting
self-avoiding polymers~\cite{FranchiniBalzan2020}.
These correspondences concern the distributions of counts; applying the
autonomous HLS theory requires control of the limiting urn function
and of the relevant boundary regimes~\cite{Franchini_HLS_review}.

For continuous urn functions, the basic asymptotic picture is
given by almost sure convergence of the color proportion to a
fixed point~\cite{HLS1980,Arthur1986}. Which limiting behaviors
are possible depends on the local shape of the response, as
illustrated by convergence to points of tangency~\cite{Pemantle1991}
and by critical behavior in nonlinear urns~\cite{NakayamaMori2021}.
Together with the characterization of the exchangeable cases
\cite{HLS1987}, these results describe typical outcomes and special
structures of reinforced growth. The fixed points alone do not
determine the exponential cost of other final compositions. Our
concern is how weighting these compositions changes the balance
between growth histories. This leads to the growth of exponential
moments and hence to the free energy studied below.

We consider two colors, black and white, and a continuous urn function
\begin{equation}
\pi:[0,1]\to[0,1].
\label{eq:binary_urn}
\end{equation}
When the urn contains $n$ balls, we denote the proportion of black
balls by $\psi_n$: the numbers of balls of the two colors are thus
$n\psi_n$ and $n(1-\psi_n)$. The new ball is black with probability
$\pi(\psi_n)$ and white with probability $1-\pi(\psi_n)$.
The urn function is evaluated before the addition and remains the
same throughout the evolution. The rule includes nonlinear reinforcement
and may admit several fixed points, that is, compositions at which
the probability of adding a black ball equals the current proportion.

The urn starts with $M\geq1$ balls and grows until it contains $T>M$.
We assume a uniform distribution over the initial compositions:
\begin{equation}
\mathbb P(\psi_M=j/M)=\frac{1}{M+1},\qquad j=0,\ldots,M.
\label{eq:initialization}
\end{equation}
We assume that $M$ is sufficiently large to avoid pathological
evolutions caused by the initial condition, that is, to keep the
outer branches considered in Theorem~\ref{thm:solution}
accessible with positive probability
\cite[Corollary~2]{Franchini_URNS}. In the limit $T\to\infty$,
$M$ remains fixed.

We denote the increment in the number of black balls by $\sigma_{n+1}$:
it equals one if the new ball is black and zero if it is white.
For $M\leq n<T$, the process is defined by
\begin{equation}
(n+1)\psi_{n+1}=n\psi_n+\sigma_{n+1},\qquad \sigma_{n+1}\in\{0,1\},
\label{eq:increments}
\end{equation}
\begin{equation}
\mathbb P(\sigma_{n+1}=1\mid\psi_M,\ldots,\psi_n)=\pi(\psi_n).
\label{eq:transitions}
\end{equation}
When $0<\pi<1$ throughout $[0,1]$, the limiting function
\eqref{eq:limiting_cgf} does not depend on the initial composition
or the finite initial number of balls
\cite[Lemma~16]{Franchini_URNS}.

For this process, the normalized cumulant generating function
(CGF) of the final black mass is
\begin{equation}
\zeta_T(\beta):=\frac1T\log\mathbb E\bigl[e^{\beta T\psi_T}\bigr],
\qquad \beta\in\mathbb R.
\label{eq:finite_cgf}
\end{equation}
The expectation includes both the initial composition and the
subsequent additions. The exponential factor weights histories
according to their final composition: positive values of $\beta$
favor more black mass, while negative values favor less. This
weighted sum over histories is a partition function, whose
normalized logarithm defines the free energy in the convention
used here. We study the limiting CGF
\begin{equation}
\zeta(\beta):=\lim_{T\to\infty}\zeta_T(\beta),
\label{eq:limiting_cgf}
\end{equation}
which describes the exponential growth of these moments as the
urn grows. In the Stanley model,
exponentially weighting the return count determines the polymer partition
function, giving this correspondence a direct thermodynamic
meaning~\cite{FranchiniBalzan2018}. Magnetization and wealth are affine
functions of the appropriate urn count, so their free energies follow
by rescaling the conjugate parameter and adding the deterministic
contribution. Varying $\beta$ probes fluctuations through the response
of the free energy. This connects the construction to large-deviation
studies of growing systems~\cite{Jack2019}.

Large deviation theory supplies the probabilistic basis for this
free energy. Bryc, Minda and Sethuraman~\cite{Bryc2009} establish
a path large deviation principle and compute the limiting CGF
for a family of affine processes related to random trees.
For urn functions in the class $\mathcal U$, recalled in
Section~\ref{subsec:model}, the large deviation principle of
Ref.~\cite{Franchini_URNS} ensures the existence of $\zeta$;
the equation~\eqref{eq:cgf_equation} and the slope conditions
selecting the solution are established in
Ref.~\cite[Eq.~(2.31) and Theorem~9]{Franchini_URNS}.
In Ref.~\cite{FranchiniBalzan2023}, this theory is applied to
increasing returns, and the equation is studied for nonlinear responses.
The broader theory also includes strong laws for nonbinary
increments~\cite{BlumBrennan1980}, under conditions on the
conditional mean and moments, and large deviations for integer
increments from zero to $K$, with $K\geq1$~\cite{Franchini_K}.
For the binary model considered here, the differential equation
and slope conditions provide the characterization from which
we construct the limiting CGF.

Although this equation is of first order, it couples the derivative
through the nonlinear urn function to an exponential of the CGF
itself. Inverting the urn function therefore does not by itself
give an integral representation. At $\beta=0$, the equation is
singular and the condition $\zeta(0)=0$ must be supplemented
by the slopes prescribed by the extreme fixed points. The slopes can
differ on the two sides of the origin, while the accessible
endpoints determine the slopes at infinity. The task is thus to
construct a globally convex solution satisfying these conditions
on both half-lines, for either orientation of the outer branches.

Theorem~\ref{thm:solution} gives an exact constructive representation
of this CGF, and hence of the asymptotic free energy, under the
assumptions of Ref.~\cite[Theorem~9]{Franchini_URNS}. These require
invertibility with Lipschitz inverses on the outer branches beyond
the fixed points, without imposing monotonicity between the extreme
fixed points. Under these assumptions, three ordinary quadratures
of the urn function and a convergent limit of compositions of
radicals yield a solution parametrized by $\zeta'(\beta)$.
Evaluation at a prescribed $\beta$ therefore requires only one
scalar inversion, while the formula itself contains no derivatives
or inverses of the urn function. The representation covers both
increasing and decreasing branches, including urn functions with
zeros and interior accessible endpoints.
It therefore supplies a common route from an admissible nonlinear
response to the free energy, without solving a new differential
equation for each model.

The radical construction is obtained from a transformation to the real
form of the chordal Loewner equation. Integration with a piecewise
constant driving function and compositions of radicals are
classical~\cite{Kennedy2009,Tran2015}; here we derive the
quadratures of the HLS transformation and identify the branch
corresponding to the CGF. To control the resulting construction,
we prove, with exact quadratures, an error bound for the reduced
variable that is uniform at the partition nodes from the initial
point up to each interior point of the branch. The estimate remains
valid at the initial singularity and does not sum the errors over
successive steps.

An exact CGF is also useful for assessing approximations to the
history law. In connection with lattice field theory descriptions
of neural activity~\cite{Bardella2024}, the HLS comparison developed
in~\cite{Franchini_HLS_review} shows that truncating the action to
pair interactions can preserve the mean activity while changing
the free energy. The present construction provides a reference for
testing such reductions for admissible nonlinear responses.
The variational balance between the action and entropy of histories
also connects this description to formulations of active inference
\cite{Friston2017,Franchini_HLS_review}.

Because the representation is expressed directly in terms of a
fixed urn function, it can also be applied conditionally when that
function is random. In the quenched setting, the response is
sampled once, from a law independent of $T$, and kept fixed
throughout growth, as in a magnetic model with a frozen random
field. Different field values give different nonlinear responses,
each with its own free energy. Urns in such fixed random environments
are investigated in~\cite{BandyopadhyayThacker2022}, while the HLS
formulation and affine examples are discussed
in~\cite{Franchini_HLS_review}. If the theorem's hypotheses,
including accessibility, hold for almost every realization, the
construction gives the corresponding conditional CGF. Each is
evaluated at a common value of $\beta$ before averaging over
disorder, since the parametrization and its inversion depend on
the realization. The construction thus provides the conditional
free energies needed to extend the affine calculations to
admissible nonlinear responses.
The distinction between conditional growth and averages over systems
also underlies the review's discussion of replica symmetry breaking.
Within such a history measure, exponential weighting changes the
relative weights of the trajectory components through their conditional
partition functions. This connects free-energy calculations to overlaps
between independently sampled histories~\cite{Franchini_HLS_review}.

In Section~\ref{sec:examples}, we illustrate the construction with
two affine urns, a non-affine urn and an urn with zeros. Their
exact solutions allow the values obtained from the representation
to be checked directly.
\section{Main result}
\label{sec:solution}

\subsection{Assumptions and fixed points}
\label{subsec:model}

We assume $\pi\in\mathcal U$, the class of Ref.~\cite{Franchini_URNS}:
there is a modulus bound $f:(0,1]\to(0,\infty)$ such that
\begin{equation}
|\pi(\alpha)-\pi(\alpha')|\leq f(|\alpha-\alpha'|)
\label{eq:modulus_bound}
\end{equation}
for sufficiently small positive distances, uniformly in
$\alpha,\alpha'\in[0,1]$, and
\begin{equation}
\lim_{\varepsilon\downarrow0}
\varepsilon\int_\varepsilon^1\frac{f(u)}{u^2}\,du=0.
\label{eq:class_U}
\end{equation}
In particular, all Hölder continuous urn functions belong to $\mathcal U$.
Define
\begin{equation}
C:=\{\alpha\in[0,1]:\pi(\alpha)=\alpha\},
\qquad c_-:=\min C,\qquad c_+:=\max C.
\label{eq:fixed_points}
\end{equation}
The case $\pi\equiv0$ gives $\zeta\equiv0$, whereas $\pi\equiv1$
gives $\zeta(\beta)=\beta$. In all other cases, the accessible
endpoints are defined by
\begin{equation}
z_-^*:=\inf\{\alpha\in[0,1]:\pi(\alpha)<1\},
\label{eq:left_endpoint}
\end{equation}
\begin{equation}
z_+^*:=\sup\{\alpha\in[0,1]:\pi(\alpha)>0\}.
\label{eq:right_endpoint}
\end{equation}
The differential equation, the limiting conditions and uniqueness
of the CGF are established in Ref.~\cite[Eq.~(2.31) and
Theorem~9]{Franchini_URNS}.

\subsection{Representation of the CGF}
\label{subsec:representation}

\begin{theorem}[Solution of the CGF equation]
\label{thm:solution}
Let $\pi\in\mathcal U$, with the model and initialization of
Section~\ref{sec:introduction}. Excluding the two constant cases
already mentioned, suppose that $\pi$ is injective on each nonempty interval among
\begin{equation}
[z_-^*,c_-),\qquad (c_+,z_+^*]
\label{eq:exterior_intervals}
\end{equation}
and that the inverse of each such restriction, defined on the corresponding
image, is Lipschitz. No monotonicity is required between $c_-$ and $c_+$.

The limit \eqref{eq:limiting_cgf} exists and is finite.
The CGF $\zeta$ is continuous and convex on $\mathbb R$ and is
$C^1$ on the two open half-lines. It solves the equation
\begin{equation}
\pi\bigl(\zeta'(\beta)\bigr)
=\frac{e^{\zeta(\beta)}-1}{e^\beta-1},
\qquad \beta\neq0,
\label{eq:cgf_equation}
\end{equation}
with the normalization
\begin{equation}
\zeta(0)=0
\label{eq:cgf_initial_value}
\end{equation}
and the slope conditions
\begin{equation}
\lim_{\beta\uparrow0}\zeta'(\beta)=c_-,
\qquad \lim_{\beta\downarrow0}\zeta'(\beta)=c_+,
\label{eq:cgf_slopes_origin}
\end{equation}
\begin{equation}
\lim_{\beta\to-\infty}\zeta'(\beta)=z_-^*,
\qquad \lim_{\beta\to+\infty}\zeta'(\beta)=z_+^*.
\label{eq:cgf_slopes_infinity}
\end{equation}
Within the class of functions with this regularity and convexity,
the problem \eqref{eq:cgf_equation}--\eqref{eq:cgf_slopes_infinity}
has a unique solution.

On the half-line $\beta>0$, if $c_+=1$, then
$\zeta(\beta)=\beta$. If $c_+<1$, fix
$\alpha\in(c_+,z_+^*)$ and set
\begin{equation}
(a,\epsilon)=
\begin{cases}
(1,-1),&\pi\text{ increasing on the upper branch},\\
(c_+,1),&\pi\text{ decreasing on the upper branch}.
\end{cases}
\label{eq:anchor}
\end{equation}
In the first case, $z_+^*=1$. On the closed interval $J_\alpha$
with endpoints $a$ and $\alpha$, the following integrals are oriented from $a$:
\begin{equation}
m(u):=
\frac{\pi(u)^u[1-\pi(u)]^{1-u}}{\pi(a)^a[1-\pi(a)]^{1-a}}
\exp\biggl[-\int_a^u\log\frac{\pi(t)}{1-\pi(t)}\,dt\biggr],
\label{eq:m_quadrature}
\end{equation}
\begin{equation}
U(u):=\frac{2u-1}{m(u)}-a+1-\int_a^u\frac{dt}{m(t)},
\label{eq:U_quadrature}
\end{equation}
\begin{equation}
R(u):=\frac{u(1-u)}{m(u)^2}-a(1-a)
-\int_a^u\frac{1-2t}{m(t)^2}\,dt.
\label{eq:R_quadrature}
\end{equation}
The three functions are continuous, $m>0$, and $R$ is strictly
increasing along $J_\alpha$ from $a$ to $\alpha$.

Let $\mathcal P=(\alpha_0,\ldots,\alpha_\ell)$ be a partition
ordered along this path, with $\alpha_0=a$ and
$\alpha_\ell=\alpha$. Define
\begin{equation}
|\mathcal P|:=\max_{1\leq j\leq\ell}|\alpha_j-\alpha_{j-1}|,
\qquad W_0:=a,
\label{eq:partition}
\end{equation}
and, for $j=1,\ldots,\ell$,
\begin{equation}
W_j:=U(\alpha_{j-1})+
\epsilon\sqrt{[W_{j-1}-U(\alpha_{j-1})]^2
+R(\alpha_j)-R(\alpha_{j-1})}.
\label{eq:radical_recursion}
\end{equation}
The sign $\epsilon$ remains the one in \eqref{eq:anchor} at every step.
The following limit exists and is independent of the partitions:
\begin{equation}
W(\alpha):=\lim_{|\mathcal P|\to0}W_\ell.
\label{eq:radical_limit}
\end{equation}
More precisely, setting, for $h\geq0$,
\begin{equation}
\omega_{U,J_\alpha}(h):=
\sup_{\substack{u,v\in J_\alpha\\ |u-v|\leq h}}
|U(u)-U(v)|,
\label{eq:driver_modulus}
\end{equation}
we have
\begin{equation}
\max_{0\leq j\leq\ell}|W_j-W(\alpha_j)|
\leq\omega_{U,J_\alpha}(|\mathcal P|),\qquad W(a):=a.
\label{eq:radical_error}
\end{equation}
The bound concerns the exact quadratures
\eqref{eq:m_quadrature}--\eqref{eq:R_quadrature}; under these assumptions,
$U$ need not be Lipschitz in the parameter $\alpha$.

Finally, define
\begin{equation}
v(\alpha):=\alpha+m(\alpha)[W(\alpha)-U(\alpha)].
\label{eq:v_reconstruction}
\end{equation}
The solution is given parametrically by
\begin{equation}
\beta(\alpha)=
\log\frac{v(\alpha)[1-\pi(\alpha)]}{\pi(\alpha)[1-v(\alpha)]},
\label{eq:lambda_reconstruction}
\end{equation}
\begin{equation}
\zeta\bigl(\beta(\alpha)\bigr)
=\log\frac{1-\pi(\alpha)}{1-v(\alpha)}.
\label{eq:zeta_reconstruction}
\end{equation}
The map $\alpha\mapsto\beta(\alpha)$ is an increasing homeomorphism
from $(c_+,z_+^*)$ onto $(0,\infty)$, and
\begin{equation}
\zeta'\bigl(\beta(\alpha)\bigr)=\alpha.
\label{eq:state_parameter}
\end{equation}
For a prescribed value of $\beta>0$, it therefore suffices to invert
this scalar map. On the increasing branch,
\begin{equation}
0<\pi(\alpha)<v(\alpha)<\alpha<1;
\label{eq:increasing_order}
\end{equation}
on the decreasing branch,
\begin{equation}
0<\pi(\alpha)<\alpha<v(\alpha)<1.
\label{eq:decreasing_order}
\end{equation}

The negative half-line is obtained by applying the same construction
to the reflected urn function:
\begin{equation}
\widehat\pi(\alpha):=1-\pi(1-\alpha),
\label{eq:reflected_urn}
\end{equation}
\begin{equation}
\zeta_\pi(\beta)=\beta+\zeta_{\widehat\pi}(-\beta),
\qquad \beta<0.
\label{eq:reflection}
\end{equation}
The subscript specifies the urn function. In particular, $c_-=0$
gives $\zeta(\beta)=0$ for $\beta\leq0$.
\end{theorem}

The quadratures are finite up to every interior point of the branch;
at the endpoints, one takes the limits of the CGF without requiring
$R$ to remain finite there.

\section{Proof}
\label{sec:proof}

\subsection{Branch selection and parametrization}
\label{subsec:proof_branch}
Uniqueness for the boundary value problem
\eqref{eq:cgf_equation}--\eqref{eq:cgf_slopes_infinity}
is given by Theorem~9 of Ref.~\cite{Franchini_URNS}.
Let $I_\pi$ be the rate function, with compact sublevel sets,
for the paths in Theorem~1 of Ref.~\cite{Franchini_URNS},
on the set $\mathcal Q$ of absolutely continuous paths $\varphi:[0,1]\to[0,1]$ with $\varphi(0)=0$ and
$0\leq\dot\varphi\leq1$ almost everywhere. This set is compact in the uniform topology.
We apply Corollary~2 of Ref.~\cite{Franchini_URNS} to
each initial composition with a nondegenerate accessible interval.
If this interval reduces to a single value, the limiting conditional
CGF is $\beta$ times that value, which is a fixed point of $\pi$
by the recursion for the mean.
The limiting CGF of the finite mixture
\eqref{eq:initialization} is the maximum of the conditional CGFs,
since its weights are positive and independent of $T$.
For $\beta>0$, in the variational problem with cost $I_\pi$, it suffices
to consider paths between the lines of slope $c_+$ and
$z_+^*$: raising the portions below the first line onto it makes their
cost zero and does not decrease the final value; no path of finite
cost exceeds the second line, since $\pi=0$ above it.
The accessibility assumed in Section~\ref{sec:introduction}
allows us to use the same cost $I_\pi$ on these paths.
For $\beta<0$, the argument is symmetric between the lines of slope
$z_-^*$ and $c_-$. Continuity of the endpoint evaluation map, the
contraction principle and Varadhan's lemma therefore give
\begin{equation}
\zeta(\beta)=\sup_{\varphi\in\mathcal Q}
\{\beta\varphi(1)-I_\pi[\varphi]\}.
\label{eq:cgf_path_variational}
\end{equation}
In particular, the limit exists and is finite and convex, with
$\zeta(0)=0$ and slopes in $[0,1]$.
The entropy integrand of $I_\pi$ is nonnegative and vanishes
exactly when $\dot\varphi(\tau)=\pi(\varphi(\tau)/\tau)$ almost everywhere.
Thus the path $\varphi(\tau)=c_+\tau$ has zero cost.
If $c_+=1$, \eqref{eq:cgf_path_variational} and $\psi_T\leq1$
give $\zeta(\beta)=\beta$ for $\beta>0$.
Henceforth assume $c_+<1$.

No zero-cost path can end above $c_+$.
Indeed, writing $\varphi(\tau)=\tau u(\log\tau)$ gives
\begin{equation}
u'(t)=\pi(u(t))-u(t),\qquad t<0.
\label{eq:zero_cost_log_time}
\end{equation}
If $u(0)=s>c_+$, the backward evolution makes $u$ increase at a
rate bounded away from zero on $[s,1]$, contradicting $u\leq1$.
Compactness of the sublevel sets of $I_\pi$ and
\eqref{eq:cgf_path_variational} give
\begin{equation}
\lim_{\beta\downarrow0}\zeta'(\beta)=c_+.
\label{eq:slope_at_zero}
\end{equation}
For now, slope limits are taken through points of
differentiability; the regularity proved below makes them ordinary limits.
On the open upper branch, continuity and injectivity give
$0<\pi(\alpha)<\alpha<1$. For the linear paths
$\varphi(\tau)=\alpha\tau$, the cost is finite:
\begin{equation}
I_\pi[\varphi]
=\alpha\log\frac{\alpha}{\pi(\alpha)}
+(1-\alpha)\log\frac{1-\alpha}{1-\pi(\alpha)}<\infty,
\qquad c_+<\alpha<z_+^*.
\label{eq:linear_path_cost}
\end{equation}
Moreover,
\begin{equation}
T\psi_T\leq Tz_+^*+M.
\label{eq:finite_mass_support}
\end{equation}
If $z_+^*<1$, the bound follows by induction starting from $T=M$
because $\pi=0$ above $z_+^*$; if $z_+^*=1$, it is immediate. Applying \eqref{eq:cgf_path_variational}
to these linear paths and using the upper bound yields
$\zeta(\beta)/\beta\to z_+^*$ as $\beta\to\infty$.
Convexity then implies
\begin{equation}
\lim_{\beta\to\infty}\zeta'(\beta)=z_+^*.
\label{eq:slope_at_infinity}
\end{equation}

Fix a point of differentiability $\beta>0$ and define the
probability measure $\mathbb P_{\beta,T}$ by
\begin{equation}
\frac{d\mathbb P_{\beta,T}}{d\mathbb P}
=\exp\{\beta T\psi_T-T\zeta_T(\beta)\}.
\label{eq:cgf_tilted_law}
\end{equation}
For every fixed neighborhood of $\zeta'(\beta)$, differentiability of
$\zeta$ allows us to choose $h>0$ small enough for the
Chernoff bounds at $\beta+h$ and $\beta-h$ to make
the probability of its complement under this measure exponentially small
as $T\to\infty$. Hence $\psi_T$ converges in probability to
$\zeta'(\beta)$ under $\mathbb P_{\beta,T}$.
Continuity and boundedness of $\pi$ imply
\begin{equation}
\lim_{T\to\infty}\mathbb E_{\beta,T}[\pi(\psi_T)]
=\pi(\zeta'(\beta)).
\label{eq:cgf_tilted_limit}
\end{equation}
The transition rule gives the exact recursion
\begin{equation}
(T+1)\zeta_{T+1}(\beta)-T\zeta_T(\beta)
=\log[1+(e^\beta-1)\mathbb E_{\beta,T}[\pi(\psi_T)]].
\label{eq:cgf_exact_recursion}
\end{equation}
The right-hand side converges by \eqref{eq:cgf_tilted_limit}.
Summing the recursion from $T=M$ and taking Ces\`aro averages
identifies the limit as $\zeta(\beta)$: the initial term,
bounded in absolute value by $|\beta|M$, vanishes after
normalization. This yields
\eqref{eq:cgf_equation} at every point of differentiability.

The slope bounds place $\zeta'$ in $[c_+,z_+^*]$
almost everywhere. Let $g$ be the inverse on the upper branch, extended
continuously to $c_+$; the extension remains Lipschitz.
The ratio in the following equation belongs to the closed image
of the branch: this holds at points of differentiability by
\eqref{eq:cgf_equation}, and elsewhere by continuity. Therefore
\begin{equation}
\zeta'(\beta)=g\biggl(\frac{e^{\zeta(\beta)}-1}{e^\beta-1}\biggr)
\qquad\text{for almost every }\beta>0.
\label{eq:cgf_ae_inverse}
\end{equation}
The right-hand side is continuous. Since $\zeta$ is locally
absolutely continuous, integration shows that it is $C^1$
on $(0,\infty)$ and that \eqref{eq:cgf_equation} holds everywhere.
The same argument applies to the reflected urn on the negative
half-line. If $c_-=0$, the zero path has zero cost and
\eqref{eq:cgf_path_variational} gives $\zeta(\beta)=0$ for $\beta<0$.
The upper branch is nonempty: otherwise, continuity and
$\pi(c_+)=c_+$ would imply $\pi\equiv0$. On this branch,
\begin{equation}
0<\pi(\alpha)<\alpha,\qquad c_+<\alpha<z_+^*.
\label{eq:upper_branch_sign}
\end{equation}
The second inequality follows from the absence of fixed points beyond
$c_+$; the first follows from strict monotonicity and the definition of
$z_+^*$. If $\pi$ is increasing, an endpoint $z_+^*<1$ would force
$\pi(z_+^*)=0$, which is impossible. Thus $z_+^*=1$ and
\begin{equation}
0\leq c_+<\pi(1)<1.
\label{eq:increasing_endpoint}
\end{equation}
If $\pi$ is decreasing, instead,
\begin{equation}
0\leq \pi(z_+^*)<c_+<1.
\label{eq:decreasing_endpoint}
\end{equation}
In both cases $0<\pi(a)<1$ and $0<\pi<1$ throughout $J_\alpha$.
The quadratures in the statement are therefore finite and continuous on
this closed interval, with $m>0$ and
\begin{equation}
m(a)=1,\qquad U(a)=a,\qquad R(a)=0.
\label{eq:anchor_values}
\end{equation}

Set $x(\beta):=\zeta'(\beta)$.
The $C^1$ regularity of the CGF makes $x$ continuous, while convexity makes it nondecreasing. We show that
\begin{equation}
c_+<x(\beta)<z_+^*,\qquad \beta>0.
\label{eq:strict_interior}
\end{equation}
If it reached either endpoint at a finite parameter value, it would
be constant on an open interval extending towards zero or infinity.
When the constant value is $d\in(0,1)$, the equation would imply
\begin{equation}
\zeta(\beta)=\log[1-\pi(d)+\pi(d)e^\beta]
\label{eq:constant_slope_contradiction}
\end{equation}
on that interval. For $0<\pi(d)<1$, the right-hand side is strictly
convex, whereas $\zeta$ is affine. For $\pi(d)=0$ or $1$, its
derivative would be zero or one, respectively, still different from $d$.
Reaching $z_+^*=1$ would require $\pi(1)=1$, which is ruled out
by $c_+<1$.

It remains to rule out reaching $c_+=0$.
In this case the branch is increasing. Its inverse $g$ extends
to zero with $g(0)=0$, and a Lipschitz constant $L$ for it gives
\mbox{$0\leq g(q)\leq Lq$}. If $x(\beta_0)=0$, convexity and
$\zeta(0)=0$ imply $\zeta(\beta_0)=0$. For every
$B>\beta_0$, the equation and $0\leq\zeta(\beta)\leq\beta$
give, on $[\beta_0,B]$,
\begin{equation}
0\leq\zeta'(\beta)
\leq\frac{L e^B}{e^{\beta_0}-1}\,\zeta(\beta).
\label{eq:zero_endpoint_gronwall}
\end{equation}
Gronwall's inequality forces $\zeta=0$ on every such
interval, contradicting \eqref{eq:slope_at_infinity}.
This proves \eqref{eq:strict_interior}.

Define
\begin{equation}
q(\beta):=\pi(x(\beta))
=\frac{e^{\zeta(\beta)}-1}{e^\beta-1},
\label{eq:q_parameter}
\end{equation}
\begin{equation}
v_\beta:=\frac{q(\beta)e^\beta}
{1-q(\beta)+q(\beta)e^\beta}.
\label{eq:tilted_probability}
\end{equation}
We have $0<q<v_\beta<1$, $x=g(q)$, and $q$ is $C^1$ by the
second expression in \eqref{eq:q_parameter}. Differentiating this expression gives
\begin{equation}
q'(\beta)=
\frac{1-q(\beta)+q(\beta)e^\beta}{e^\beta-1}
[x(\beta)-v_\beta].
\label{eq:q_derivative}
\end{equation}
Moreover, $q$ is nondecreasing if $\pi$ is increasing and nonincreasing
if $\pi$ is decreasing.

The derivative $q'$ cannot vanish. Indeed, if
$q'(\beta_0)=0$, the Lipschitz property of $g$ implies
$x(\beta_0+h)-x(\beta_0)=o(h)$.
Differentiating \eqref{eq:tilted_probability} at that point gives
\begin{equation}
\left.\frac{d}{d\beta}[x(\beta)-v_\beta]
\right|_{\beta=\beta_0}
=-v_{\beta_0}(1-v_{\beta_0})<0.
\label{eq:nonvanishing_derivative}
\end{equation}
Since \eqref{eq:q_derivative} would imply
$x(\beta_0)=v_{\beta_0}$, the sign of $q'$ would change from
positive to negative across $\beta_0$, contradicting
the monotonicity of $q$. Since $q'$ does not vanish,
$x=g(q)$ is also strictly increasing in $\beta$.
This proves the required homeomorphism property for the parameter
$\alpha=x(\beta)$ and, by \eqref{eq:q_derivative}, the orderings
\eqref{eq:increasing_order} and \eqref{eq:decreasing_order}.

Express the probability \eqref{eq:tilted_probability} as a function of
the state by $v(x(\beta))=v_\beta$. We will show below that this
function satisfies the reconstruction formula \eqref{eq:v_reconstruction}.
Set
\begin{equation}
A(q):=\frac{g(q)-q}{q(1-q)}.
\label{eq:A_coefficient}
\end{equation}
This function is continuous and positive in the interior.
The logarithmic derivative of \eqref{eq:tilted_probability} gives
\begin{equation}
\frac{dv_\beta}{d\beta}
=v_\beta(1-v_\beta)\biggl[1+\frac{q'(\beta)}{q(\beta)[1-q(\beta)]}\biggr].
\label{eq:v_lambda_derivative}
\end{equation}
The definition of $v_\beta$ also gives
\begin{equation}
\frac{1-q(\beta)+q(\beta)e^\beta}{e^\beta-1}
=\frac{q(\beta)[1-q(\beta)]}{v_\beta-q(\beta)}.
\label{eq:prefactor_identity}
\end{equation}
Dividing \eqref{eq:v_lambda_derivative} by the nonzero derivative in
\eqref{eq:q_derivative}, we obtain
\begin{equation}
\frac{d}{dq}v(g(q))
=A(q)\frac{v(g(q))[1-v(g(q))]}{g(q)-v(g(q))}.
\label{eq:rational_ode}
\end{equation}

\subsection{Quadratures and reduction}
\label{subsec:proof_reduction}
To differentiate the quadratures without differentiating $\pi$,
compose the functions $m,U,R$ in the statement with $g$.
They are absolutely continuous on every closed interval between the initial
point $q_a:=\pi(a)$ and an interior point. Indeed, $g$ is Lipschitz and
the primitives of the integrands in
\eqref{eq:m_quadrature}--\eqref{eq:R_quadrature} are $C^1$ on
$J_\alpha$. The explicit terms are also absolutely continuous after
the substitution $\pi(g(q))=q$.
In the logarithmic derivative of \eqref{eq:m_quadrature}, the terms
containing $g'$ cancel, leaving, almost everywhere,
\begin{equation}
\frac{d}{dq}\log m(g(q))
=\frac{g(q)}q-\frac{1-g(q)}{1-q}=A(q).
\label{eq:log_m_derivative}
\end{equation}
Thus
\begin{equation}
\frac{d}{dq}m(g(q))=A(q)m(g(q)).
\label{eq:m_derivative}
\end{equation}
Similarly,
\begin{equation}
\frac{d}{dq}U(g(q))
=\frac{g'(q)-A(q)[2g(q)-1]}{m(g(q))}
\label{eq:U_derivative}
\end{equation}
almost everywhere, and
\begin{equation}
\frac{d}{dq}R(g(q))
=-\frac{2A(q)g(q)[1-g(q)]}{m(g(q))^2}.
\label{eq:R_derivative}
\end{equation}
The identities \eqref{eq:m_derivative} and \eqref{eq:R_derivative}
initially hold almost everywhere; their continuous right-hand sides
then show that both compositions are $C^1$.
Absolute continuity of $\pi$ is not required.

Introduce
\begin{equation}
w(q):=U(g(q))+\frac{v(g(q))-g(q)}{m(g(q))},
\label{eq:w_transform}
\end{equation}
\begin{equation}
c(q):=\frac{A(q)g(q)[1-g(q)]}{m(g(q))^2}.
\label{eq:c_coefficient}
\end{equation}
Using the algebraic identity
\begin{equation}
\frac{v(1-v)}{g-v}=v+g-1+\frac{g(1-g)}{g-v},
\label{eq:polynomial_division}
\end{equation}
the derivatives \eqref{eq:rational_ode}--\eqref{eq:U_derivative} give
\begin{equation}
\frac{dw}{dq}=-\frac{c(q)}{w(q)-U(g(q))}.
\label{eq:w_ode}
\end{equation}
The identity holds almost everywhere; the denominator
does not vanish in the interior, by
\eqref{eq:increasing_order}--\eqref{eq:decreasing_order}, and the
continuous right-hand side shows that $w$ is $C^1$.

On the increasing branch we move downwards in $\alpha$ from $a=1$,
that is, from infinity towards zero in $\beta$; on the decreasing
branch we move upwards in $\alpha$ from $a=c_+$, that is, from zero
towards infinity in $\beta$. Along both paths, $q$ decreases from
$q_a$. Since $c>0$ in the interior, \eqref{eq:R_derivative} proves
that $R$ is strictly increasing along the path.
At the starting endpoints,
\begin{equation}
\lim_{\beta\to\infty}v_\beta=1
\qquad\text{on the increasing branch},
\label{eq:increasing_initial_v}
\end{equation}
\begin{equation}
\lim_{\beta\downarrow0}v_\beta=c_+
\qquad\text{on the decreasing branch}.
\label{eq:decreasing_initial_v}
\end{equation}
These formulas follow from \eqref{eq:tilted_probability} and the
limits already established for $q$. In both cases, $w$ extends
continuously to the initial point with value $a$.

We now use $t=R(u)$, with $u\in J_\alpha$, as time. Since $R$ is
continuous and strictly increasing along $J_\alpha$, the identities
\begin{equation}
y(R(u)):=w(\pi(u)),\qquad \xi(R(u)):=U(u),
\label{eq:time_functions}
\end{equation}
define $y$ and $\xi$ uniquely on $[0,R(\alpha)]$.
Differentiating with respect to $q$, with $t=R(g(q))$,
equations \eqref{eq:R_derivative} and \eqref{eq:w_ode} give
\begin{equation}
y'(t)=\frac{1}{2[y(t)-\xi(t)]},\qquad 0<t<R(\alpha),
\label{eq:normal_ode}
\end{equation}
with
\begin{equation}
y(0)=\xi(0)=a.
\label{eq:singular_initial_value}
\end{equation}
For $t>0$, we have $y<\xi$ on the increasing branch and $y>\xi$
on the decreasing branch. The time is finite for every interior
$\alpha$, and the denominator does not vanish before that time
is reached. The solution on the prescribed side of the driving function is obtained from the CGF.
After rescaling time, \eqref{eq:normal_ode} is the chordal Loewner
equation in a real variable \cite{Kennedy2009}. Integrating with
$\xi$ constant on each step gives the radical used in the construction.

\subsection{Comparison with a singular initial condition}
\label{subsec:proof_comparison}
We first treat $y<\xi$. For $h\geq0$, define
\begin{equation}
P_{b,h}(z):=b-\sqrt{(b-z)_+^2+h},
\qquad r_+:=\max\{r,0\},
\label{eq:projected_flow}
\end{equation}
\begin{equation}
N_{b,h}(Z):=b-\sqrt{(Z-b)^2+h}.
\label{eq:fixed_sign_flow}
\end{equation}
The first map is nondecreasing in $z$ and satisfies
\begin{equation}
P_{b+d,h}(z+d)=P_{b,h}(z)+d.
\label{eq:translation_covariance}
\end{equation}
For $h>0$, it is the solution with constant driving function $b$
starting from $\min\{z,b\}$, including when the initial condition
equals $b$. The recursion uses $N_{b,h}$.

Consider a step of duration $h>0$, shifting its time origin to zero.
Assume
\begin{equation}
y_0:=y(0)\leq b,\qquad |\xi(t)-b|\leq\delta,
\qquad 0\leq t\leq h.
\label{eq:one_step_assumptions}
\end{equation}
Then
\begin{equation}
P_{b-\delta,h}(y_0)\leq y(h)\leq P_{b+\delta,h}(y_0).
\label{eq:exact_flow_comparison}
\end{equation}
To verify the bound also for a singular initial condition,
set $B=b-\delta$, $D=b+\delta$, and consider the functions
$l(t)=P_{B,t}(y_0)$ and $u(t)=P_{D,t}(y_0)$.
They are continuous at zero and smooth for $t>0$, with
$l(0)=\min\{y_0,B\}\leq y_0=u(0)$. Write
\begin{equation}
F_b(z):=\frac{1}{2(z-b)},\qquad z<b.
\label{eq:comparison_field}
\end{equation}
The field $F_b$ is decreasing in $z$ and increasing in $b$
wherever it is defined. If $l(t)>y(t)$, then
$y(t)<l(t)<B\leq\xi(t)$, hence
\begin{equation}
(l-y)'=F_B(l)-F_{\xi(t)}(y)
\leq F_B(l)-F_B(y)<0.
\label{eq:lower_barrier_derivative}
\end{equation}
If instead $y(t)>u(t)$, then
$u(t)<y(t)<\xi(t)\leq D$ and
\begin{equation}
(y-u)'=F_{\xi(t)}(y)-F_D(u)
\leq F_D(y)-F_D(u)<0.
\label{eq:upper_barrier_derivative}
\end{equation}
On every open component of the set where either difference is
positive, that difference would therefore have a nonpositive initial
value and a negative derivative, which is impossible. Continuity at
zero allows the same argument for a component starting at zero,
without requiring the solution to be differentiable there. This proves
\eqref{eq:exact_flow_comparison}.

Now let $Z$ be an approximate initial value with $|Z-y_0|\leq\delta$.
The triangle inequality and $y_0\leq b$ give
\begin{equation}
(b-y_0-\delta)_+\leq|Z-b|\leq b-y_0+\delta.
\label{eq:distance_bounds}
\end{equation}
Therefore
\begin{equation}
P_{b,h}(y_0-\delta)\leq N_{b,h}(Z)
\leq P_{b,h}(y_0+\delta).
\label{eq:numerical_flow_comparison}
\end{equation}
Identity \eqref{eq:translation_covariance} gives
\begin{equation}
P_{b,h}(y_0-\delta)=P_{b+\delta,h}(y_0)-\delta,
\label{eq:lower_translation}
\end{equation}
\begin{equation}
P_{b,h}(y_0+\delta)=P_{b-\delta,h}(y_0)+\delta.
\label{eq:upper_translation}
\end{equation}
Combining \eqref{eq:exact_flow_comparison} with
\eqref{eq:numerical_flow_comparison}, we obtain
\begin{equation}
|N_{b,h}(Z)-y(h)|\leq\delta.
\label{eq:one_step_error}
\end{equation}
Thus the same error bound passes from one step to the next,
even when $Z>b$.

\subsection{Convergence and error control}
\label{subsec:proof_convergence}
For the partition in the statement, set
\begin{equation}
t_j:=R(\alpha_j),\qquad
\delta:=\omega_{U,J_\alpha}(|\mathcal P|).
\label{eq:discrete_times}
\end{equation}
During the step from $t_{j-1}$ to $t_j$, the unique $u\in J_\alpha$
satisfying $R(u)=t$ lies between $\alpha_{j-1}$ and $\alpha_j$.
The definition of $\xi$ therefore gives
\begin{equation}
|\xi(t)-U(\alpha_{j-1})|\leq\delta,
\qquad t_{j-1}\leq t\leq t_j.
\label{eq:driver_step_bound}
\end{equation}
On the increasing branch, $\epsilon=-1$, and the recursion
\eqref{eq:radical_recursion} is precisely
$N_{b,h}$ with $b=U(\alpha_{j-1})$ and $h=t_j-t_{j-1}$.
Since $W_0=y(0)$, the bound \eqref{eq:one_step_error} applies
inductively to every step and gives
\begin{equation}
\max_{0\leq j\leq\ell}|W_j-y(t_j)|\leq\delta.
\label{eq:inductive_error}
\end{equation}
On the decreasing branch, the substitution
$y\mapsto-y$, $\xi\mapsto-\xi$ preserves
\eqref{eq:normal_ode} and reduces the problem to the case just proved.
Returning to the original sign gives exactly the positive square root
prescribed by $\epsilon=1$, with the same bound.

Uniform continuity of $U$ on $J_\alpha$ implies
$\delta\to0$ as $|\mathcal P|\to0$, without any condition
on the ratios of the lengths of different steps.
Hence the limit exists, is independent of the partitions, and is
\begin{equation}
W(\alpha)=w(\pi(\alpha))=y(R(\alpha)).
\label{eq:limit_identification}
\end{equation}
The same identification at each node turns
\eqref{eq:inductive_error} into the bound
\eqref{eq:radical_error} in the statement.

The same comparison applies to any other continuous solution of
\eqref{eq:normal_ode} with the same initial value and lying strictly
on the same side of $\xi$ for $t>0$. It satisfies
\eqref{eq:inductive_error} relative to the same radicals and therefore
coincides with \eqref{eq:limit_identification} as the mesh size tends
to zero. The reduced solution is thus unique on this branch.

\subsection{Reconstruction and reflection}
\label{subsec:proof_reconstruction}
The identities \eqref{eq:w_transform} and \eqref{eq:limit_identification}
recover \eqref{eq:v_reconstruction}. Solving
\eqref{eq:tilted_probability} for the exponential and using
\eqref{eq:cgf_equation}, we obtain, respectively,
\begin{equation}
e^\beta=\frac{v_\beta[1-q(\beta)]}{q(\beta)(1-v_\beta)},
\label{eq:exponential_lambda}
\end{equation}
\begin{equation}
e^{\zeta(\beta)}=\frac{1-q(\beta)}{1-v_\beta}.
\label{eq:exponential_zeta}
\end{equation}
With $\alpha=x(\beta)$, $q(\beta)=\pi(\alpha)$ and
$v_\beta=v(\alpha)$, these are precisely
\eqref{eq:lambda_reconstruction} and \eqref{eq:zeta_reconstruction}.
The arguments of their logarithms are positive by the orderings
already proved; $\beta(\alpha)>0$ follows from $v(\alpha)>\pi(\alpha)$.
Strict monotonicity and the limits
\eqref{eq:slope_at_zero}--\eqref{eq:slope_at_infinity} complete
the parametrization of the entire positive half-line.

Finally, the proportion of white mass $1-\psi_n$ has urn function
$\widehat\pi$ and the same discrete uniform initialization
\eqref{eq:initialization}. For every finite $T$, we already have
\begin{equation}
\zeta_{\pi,T}(\beta)
=\beta+\zeta_{\widehat\pi,T}(-\beta).
\label{eq:finite_reflection}
\end{equation}
Passing to the limit gives \eqref{eq:reflection}. Reflection
preserves the class $\mathcal U$ and the Lipschitz property of the
inverses, interchanging the two outer branches. If $c_-=0$, the
largest fixed point of the reflected function is one, and the case
already treated gives $\zeta(\beta)=0$ for $\beta<0$.
The normalization $\zeta(0)=0$ completes the solution on
$\mathbb R$.
\qed

\section{Examples and computation of the CGF}
\label{sec:examples}

We compare the recursion \eqref{eq:radical_recursion}
with four explicit CGFs.

\subsection{Two affine urns}
\label{subsec:affine}

For comparison, we use two exact affine solutions from
Corollary~12 of Ref.~\cite{Franchini_URNS}, as corrected in the
accompanying corrigendum. We first consider
\begin{equation}
\pi(\alpha)=\frac{1+\alpha}{3}.
\label{eq:affine_increasing}
\end{equation}
The fixed point is $c_-=c_+=1/2$, the inverse is Lipschitz with
constant $3$, and the positive branch starts at $a=1$ with
$\epsilon=-1$. The three quadratures become
\begin{equation}
m(\alpha)=\frac{2}{(1+\alpha)(2-\alpha)},
\label{eq:affine_m}
\end{equation}
\begin{equation}
U(\alpha)=-\frac{5\alpha^3}{6}+\frac{5\alpha^2}{4}
+\frac{\alpha}{2}+\frac1{12},
\label{eq:affine_U}
\end{equation}
\begin{equation}
R(\alpha)=\frac{\alpha^2}{2}-\frac{5\alpha^3}{6}
+\frac{\alpha^5}{2}-\frac{\alpha^6}{6}.
\label{eq:affine_R}
\end{equation}
The limit of the radicals is
\begin{equation}
W(\alpha)=-\frac{\alpha^3}{3}+\alpha^2+\frac13,
\qquad 1/2<\alpha<1.
\label{eq:affine_W}
\end{equation}
Substitution into \eqref{eq:v_reconstruction} gives
\begin{equation}
v(\alpha)=\frac{1+\alpha}{2(2-\alpha)},
\label{eq:affine_v}
\end{equation}
and reconstruction yields
\begin{equation}
\zeta(\beta)=\beta-\log\frac{3-e^{-\beta}}2,
\qquad \beta\geq0.
\label{eq:affine_positive_cgf}
\end{equation}
The function in \eqref{eq:affine_W}, expressed in the time variable
$t=R(\alpha)$, satisfies \eqref{eq:normal_ode}, has $W(1)=1$
and satisfies $W<U$ for $1/2<\alpha<1$. By the uniqueness proved in
Section~\ref{subsec:proof_convergence}, this function coincides
with the limit of the radicals; \eqref{eq:affine_positive_cgf} is therefore
the CGF. Since $\pi(1-\alpha)=1-\pi(\alpha)$, reflection completes the solution:
\begin{equation}
\zeta(\beta)=-\log\frac{3-e^\beta}2,
\qquad \beta\leq0.
\label{eq:affine_negative_cgf}
\end{equation}

For the second example, take
\begin{equation}
\pi(\alpha)=\frac34-\frac\alpha2.
\label{eq:affine_decreasing}
\end{equation}
The fixed point is again $1/2$, but the inverse has Lipschitz constant $2$
and the positive branch starts at $a=1/2$, with $\epsilon=1$.
The first quadrature is
\begin{equation}
m(\alpha)=\{4\pi(\alpha)[1-\pi(\alpha)]\}^{3/2}.
\label{eq:decreasing_m}
\end{equation}
The exact solution is now
\begin{equation}
\zeta(\beta)=\frac\beta2+2\log\cosh\frac\beta4,
\qquad \beta\in\mathbb R.
\label{eq:decreasing_cgf}
\end{equation}
Indeed, setting $t=e^{\beta/2}$ gives
\begin{equation}
\zeta'(\beta)=\frac{t}{1+t},
\qquad e^{\zeta(\beta)}=\frac{(1+t)^2}{4},
\label{eq:decreasing_check}
\end{equation}
and equation \eqref{eq:cgf_equation} follows by substitution.
For $\beta>0$, the derivative is strictly increasing and
\eqref{eq:q_derivative} gives $v>\alpha$; moreover, $\alpha$ and $v$
both tend to $1/2$ as $\beta\downarrow0$. The transformed
solution therefore starts at the prescribed initial point and remains above
$\xi$. Uniqueness from Section~\ref{subsec:proof_convergence}
identifies the CGF on the positive half-line; the symmetry of the urn
function gives the negative half-line by reflection.
Figure~\ref{fig:affine} compares the two exact CGFs
with the values obtained from the recursion.

\afterpage{%
\newgeometry{left=20mm,right=20mm,top=18mm,bottom=18mm}
\begingroup
\singlespacing
\setlength{\abovecaptionskip}{10pt}
\setlength{\belowcaptionskip}{0pt}
\begin{figure}[p]
\centering
\includegraphics[width=\linewidth,height=0.72\textheight,keepaspectratio]{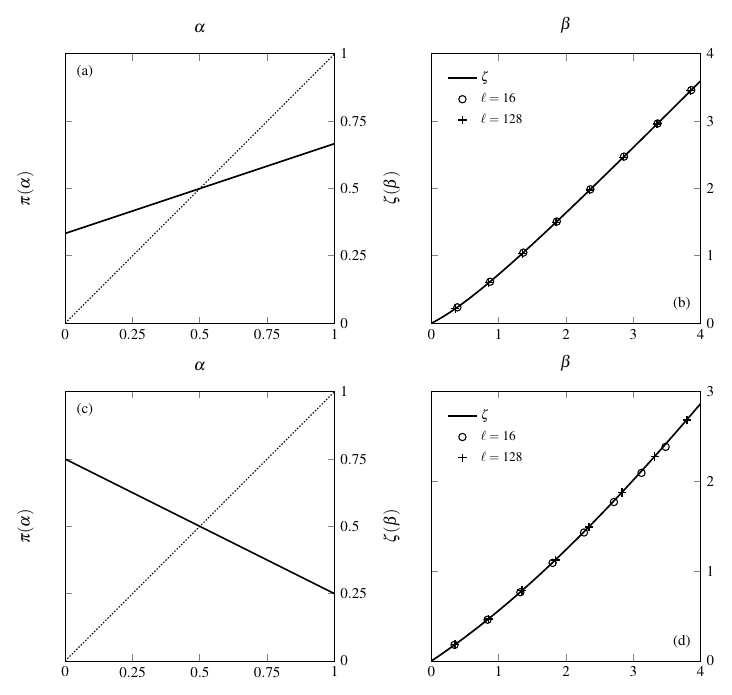}
\caption{Two affine urns. Top: the increasing function
\eqref{eq:affine_increasing}; bottom: the decreasing function
\eqref{eq:affine_decreasing}. On the left, the urn function
(solid line) and the diagonal (dotted line) intersect at the fixed
point. On the right, the exact CGF is compared with
the points obtained from \eqref{eq:radical_recursion}
using $\ell=16$ intervals (circles) and $\ell=128$ (crosses), uniformly spaced in
$\alpha$ between the initial point and the evaluation point. The quadratures are
computed on a finer grid. The upper example starts
at the boundary $a=1$ with the negative square root; the lower example
starts at the fixed point $a=1/2$ with the positive square root.}
\label{fig:affine}
\end{figure}
\clearpage
\endgroup
\restoretextgeometry
}

\subsection{A non-affine urn}
\label{subsec:nonlinear}

Consider the increasing function
\begin{equation}
\pi(\alpha)=\frac{\sqrt\alpha+\sqrt{1-\alpha}}
{\sqrt2+2\sqrt{1-\alpha}},\qquad 0\leq\alpha\leq1.
\label{eq:nonlinear_urn}
\end{equation}
It is Hölder continuous with exponent $1/2$, hence belongs to $\mathcal U$, and
satisfies \mbox{$\pi(1-\alpha)=1-\pi(\alpha)$}. Its inverse is
\begin{equation}
g(q)=\frac{(-2q^2+4q-1)^2}{2[q^2+(1-q)^2]^2},
\qquad 1-1/\sqrt2\leq q\leq1/\sqrt2.
\label{eq:nonlinear_inverse}
\end{equation}
With $u=2q-1$, we obtain
\begin{equation}
g'(q)=\frac{4(1-6u^2+u^4)}{(1+u^2)^3},
\qquad 0\leq g'(q)\leq4.
\label{eq:nonlinear_inverse_derivative}
\end{equation}
The inverse is therefore Lipschitz, although $\pi$ has no finite derivative
at the endpoints. Moreover,
\begin{equation}
g(q)-q=\frac{u(3-6u^2-u^4)}{2(1+u^2)^2}.
\label{eq:nonlinear_fixed_point}
\end{equation}
The factor in parentheses is positive throughout the interval under consideration:
the unique fixed point is $1/2$. For the positive half-line, we therefore use
$a=1$ and $\epsilon=-1$.

The CGF is
\begin{equation}
\zeta(\beta)=\frac12\log\frac{1+e^{2\beta}}2,
\qquad \beta\in\mathbb R.
\label{eq:nonlinear_cgf}
\end{equation}
To verify this, set $t=e^\beta$; then
\begin{equation}
\alpha=\zeta'(\beta)=\frac{t^2}{1+t^2}.
\label{eq:nonlinear_slope}
\end{equation}
Rationalizing the quotient in the HLS equation gives
\begin{equation}
\frac{e^{\zeta(\beta)}-1}{e^\beta-1}
=\frac{t+1}{\sqrt2\sqrt{1+t^2}+2}=\pi(\alpha),
\qquad \beta\neq0.
\label{eq:nonlinear_check}
\end{equation}
On the positive half-line, $q=\pi(\alpha)$ has positive derivative,
so \eqref{eq:q_derivative} gives $v<\alpha$; moreover,
$\alpha$ and $v$ both tend to one as $\beta\to\infty$.
The transformed solution remains below $\xi$ and has the prescribed initial point.
Uniqueness from Section~\ref{subsec:proof_convergence} therefore identifies
this formula as the CGF; symmetry completes the negative half-line.

\subsection{Zeros and an interior accessible endpoint}
\label{subsec:boundary}

Take
\begin{equation}
\pi(\alpha)=\max\{1-2\alpha,0\}.
\label{eq:boundary_urn}
\end{equation}
In this case
\begin{equation}
c_-=c_+=1/3,\qquad z_-^*=0,\qquad z_+^*=1/2.
\label{eq:boundary_endpoints}
\end{equation}
The inverse on the two outer branches is $g(q)=(1-q)/2$ and is Lipschitz;
injectivity in \eqref{eq:exterior_intervals} is required
only up to the endpoint $z_+^*=1/2$. On the positive half-line, the construction
starts at $a=1/3$ with $\epsilon=1$ and remains in the interval
$(1/3,1/2)$. The first quadrature becomes
\begin{equation}
m(\alpha)=3\alpha\sqrt{3(1-2\alpha)}.
\label{eq:boundary_m}
\end{equation}

The solution on the whole real line can be written as
\begin{equation}
e^{-\zeta(\beta)}=F(e^\beta-1),
\label{eq:boundary_cgf_integral}
\end{equation}
\begin{equation}
F(z):=\int_0^1\frac{du}{1+zu^2},\qquad z>-1.
\label{eq:boundary_F}
\end{equation}
In particular,
\begin{equation}
e^{-\zeta(\beta)}=
\frac{\arctan\sqrt{e^\beta-1}}{\sqrt{e^\beta-1}},
\qquad \beta>0,
\label{eq:boundary_positive_cgf}
\end{equation}
\begin{equation}
e^{-\zeta(\beta)}=
\frac{\operatorname{artanh}\sqrt{1-e^\beta}}{\sqrt{1-e^\beta}},
\qquad \beta<0,
\label{eq:boundary_negative_cgf}
\end{equation}
with $\zeta(0)=0$. Indeed, on the branches where
$\pi(\alpha)=1-2\alpha$, the HLS equation becomes
\begin{equation}
2zF'(z)+F(z)=\frac1{1+z}.
\label{eq:boundary_linear_ode}
\end{equation}
The integral \eqref{eq:boundary_F} satisfies this identity.
On each half-line, the difference between two solutions is a constant divided
by the square root of $|z|$; continuity at zero eliminates
both constants. Moreover, $F'<0$ and
\begin{equation}
\zeta'(\beta)=
\frac{e^\beta F(e^\beta-1)-1}
{2(e^\beta-1)F(e^\beta-1)},\qquad \beta\neq0,
\label{eq:boundary_slope}
\end{equation}
which gives $0<\zeta'(\beta)<1/2$: the derivative remains in the interval
where the urn function is affine. We recover
\begin{equation}
\zeta'(0)=1/3,\qquad
\lim_{\beta\to\infty}\zeta'(\beta)=1/2.
\label{eq:boundary_slope_limits}
\end{equation}
At the interior endpoint, we also have
\begin{equation}
\lim_{\beta\to\infty}[\zeta(\beta)-\beta/2]
=\log(2/\pi).
\label{eq:boundary_asymptote}
\end{equation}
In equations \eqref{eq:boundary_asymptote} and
\eqref{eq:quadrature_substitution}, $\pi$ is the circle constant.
The constant term is finite even though $\pi(z_+^*)=0$.

Figure~\ref{fig:nonlinear_boundary} compares the exact
solutions with the values from the recursion for the non-affine urn and the urn
with zeros.
In the non-affine example, the urn function remains positive and the slope of the
CGF tends to one; for the urn with zeros, the range of slopes
ends at the interior endpoint $1/2$.

\afterpage{%
\newgeometry{left=20mm,right=20mm,top=18mm,bottom=18mm}
\begingroup
\singlespacing
\setlength{\abovecaptionskip}{10pt}
\setlength{\belowcaptionskip}{0pt}
\begin{figure}[p]
\centering
\includegraphics[width=\linewidth,height=0.72\textheight,keepaspectratio]{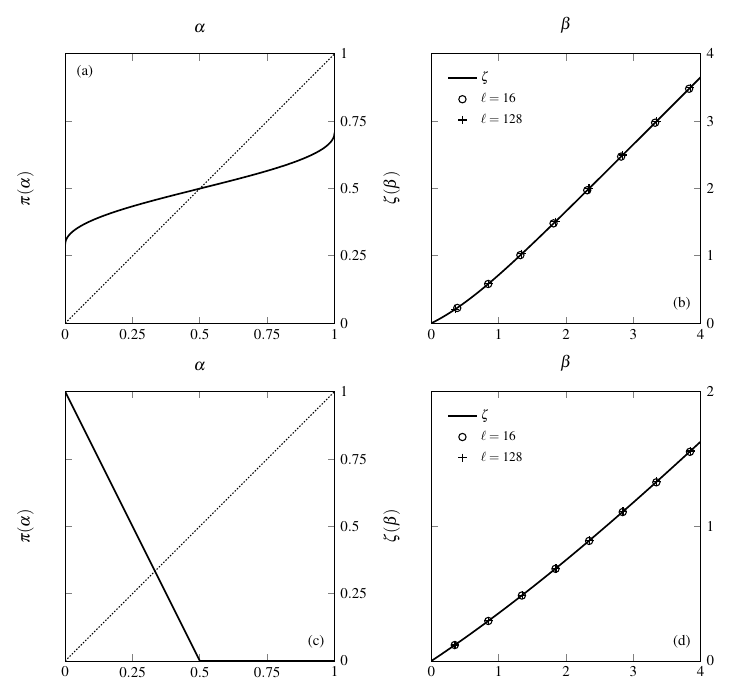}
\caption{A non-affine urn and an urn with zeros. Top:
\eqref{eq:nonlinear_urn}, an increasing function with a Lipschitz
inverse and fixed point $1/2$; the CGF
\eqref{eq:nonlinear_cgf} has limiting slope one.
Bottom: \eqref{eq:boundary_urn} vanishes
from $1/2$ onwards; the fixed point is $1/3$ and the limiting CGF
slope is $1/2$. Left: urn function and diagonal;
right: exact CGF and points computed with $\ell=16$ intervals
(circles) and $\ell=128$ (crosses), using the same conventions as in
Figure~\ref{fig:affine}. The upper example starts at $a=1$ with
the negative square root, the lower example at $a=1/3$ with the positive square root.}
\label{fig:nonlinear_boundary}
\end{figure}
\clearpage
\endgroup
\restoretextgeometry
}

\paragraph{Numerical evaluation.}
Each approximate point in the figures is obtained by substituting $W_\ell$ for $W$
in formulas \eqref{eq:v_reconstruction}--\eqref{eq:zeta_reconstruction}.
We denote the coordinates by $\beta_{\mathcal P}$ and $\zeta_{\mathcal P}$.
The vertical error relative to the exact curve is
\begin{equation}
|\zeta_{\mathcal P}(\alpha)
-\zeta(\beta_{\mathcal P}(\alpha))|.
\label{eq:numerical_graph_error}
\end{equation}
Comparison at a prescribed value of $\beta$ also requires
inverting the parametrization. For the four points corresponding
to the exact values $\beta=0.25,1,2,4$, with $\ell=2048$ intervals
uniformly spaced in $\alpha$, the maximum of
\eqref{eq:numerical_graph_error} is less than
$3.3\cdot10^{-4}$ in all examples.
The quadratures were computed using the trapezoidal rule
after the substitution
\begin{equation}
u=a+(\alpha-a)\sin^2(\pi s/2),\qquad 0\leq s\leq1.
\label{eq:quadrature_substitution}
\end{equation}
In these checks, doubling the grid from $65536$ to $131072$ intervals
changes the reconstructed coordinates by less than $8\cdot10^{-9}$.
This comparison measures the numerical sensitivity of the quadratures;
the theoretical bound \eqref{eq:radical_error} concerns exact
quadratures. On an interior interval, the reconstruction
formulas are regular functions of the admissible data. Notice that near the endpoints the logarithms may amplify the error and require a separate check.

\section*{Declaration of Generative AI}
Research, writing, and proofreading of the present manuscript were performed with AI assistance (GPT-6 Astra) under strict human guidance and scrutiny. We independently checked all mathematical arguments and take full responsibility for the content of the publication.


\begin{thebibliography}{22}
\bibitem{HLS1980}
B.~M.~Hill, D.~Lane, W.~Sudderth,
\textit{A Strong Law for Some Generalized Urn Processes},
Ann. Probab. 8(2), 214--226 (1980).
\url{https://doi.org/10.1214/aop/1176994772}.

\bibitem{Pemantle2007}
R.~Pemantle,
\textit{A survey of random processes with reinforcement},
Probab. Surv. 4, 1--79 (2007).
\url{https://doi.org/10.1214/07-PS094}.

\bibitem{Franchini_HLS_review}
S.~Franchini,
\textit{The Urn of Hill, Lane and Sudderth},
\href{https://arxiv.org/abs/2506.20826}{arXiv:2506.20826} (2026).
\url{https://doi.org/10.48550/arXiv.2506.20826}.

\bibitem{Bryc2009}
W.~Bryc, D.~Minda, S.~Sethuraman,
\textit{Large deviations for the leaves in some random trees},
Adv. Appl. Probab. 41(3), 845--873 (2009).
\url{https://doi.org/10.1239/aap/1253281066}.

\bibitem{BaurBertoin2016}
E.~Baur, J.~Bertoin,
\textit{Elephant random walks and their connection to P\'olya-type urns},
Phys. Rev. E 94, 052134 (2016).
\url{https://doi.org/10.1103/PhysRevE.94.052134}.

\bibitem{PodderRoy2026}
M.~Podder, A.~Roy,
\textit{Elephant Random Walks with Multiple Extractions and General Reinforcement Functions},
J. Theor. Probab. 39(1), 17 (2026).
\url{https://doi.org/10.1007/s10959-025-01471-4}.

\bibitem{DosiErmolievKaniovski1994}
G.~Dosi, Y.~Ermoliev, Y.~Kaniovski,
\textit{Generalized urn schemes and technological dynamics},
J. Math. Econ. 23(1), 1--19 (1994).
\url{https://doi.org/10.1016/0304-4068(94)90032-9}.

\bibitem{FranchiniBalzan2023}
S.~Franchini, R.~Balzan,
\textit{Large-deviation theory of increasing returns},
Phys. Rev. E 107, 064142 (2023).
\url{https://doi.org/10.1103/PhysRevE.107.064142}.

\bibitem{GottfriedGrosskinsky2025}
T.~Gottfried, S.~Grosskinsky,
\textit{Wages and capital returns in a generalized P\'olya urn},
J. Econ. Interact. Coord. 20, 477--518 (2025).
\url{https://doi.org/10.1007/s11403-024-00437-9}.

\bibitem{KlymkoGarrahanWhitelam2017}
K.~Klymko, J.~P.~Garrahan, S.~Whitelam,
\textit{Similarity of ensembles of trajectories of reversible and irreversible growth processes},
Phys. Rev. E 96, 042126 (2017).
\url{https://doi.org/10.1103/PhysRevE.96.042126}.

\bibitem{KhaninKhanin2001}
K.~Khanin, R.~Khanin,
\textit{A probabilistic model for the establishment of neuron polarity},
J. Math. Biol. 42, 26--40 (2001).
\url{https://doi.org/10.1007/PL00000071}.

\bibitem{Barabas2009}
B.~Barab\'as, L.~Caglioti, K.~Micskei, G.~P\'alyi,
\textit{Data-Based Stochastic Approach to Absolute Asymmetric Synthesis by Autocatalysis},
Bull. Chem. Soc. Jpn. 82(11), 1372--1376 (2009).
\url{https://doi.org/10.1246/bcsj.82.1372}.

\bibitem{BriendDevroyeLugosi2025}
S.~Briend, L.~Devroye, G.~Lugosi,
\textit{Broadcasting in random recursive dags},
ESAIM Probab. Stat. 29, 184--203 (2025).
\url{https://doi.org/10.1051/ps/2025002}.

\bibitem{Bellina2023}
A.~Bellina, C.~Castellano, P.~Pineau, G.~Iannelli, G.~De Marzo,
\textit{Effect of collaborative-filtering-based recommendation algorithms on opinion polarization},
Phys. Rev. E 108, 054304 (2023).
\url{https://doi.org/10.1103/PhysRevE.108.054304}.

\bibitem{Tria2014}
F.~Tria, V.~Loreto, V.~D.~P.~Servedio, S.~H.~Strogatz,
\textit{The dynamics of correlated novelties},
Sci. Rep. 4, 5890 (2014).
\url{https://doi.org/10.1038/srep05890}.

\bibitem{FranchiniBalzan2018}
S.~Franchini, R.~Balzan,
\textit{Random polymers and generalized urn processes},
Phys. Rev. E 98, 042502 (2018).
\url{https://doi.org/10.1103/PhysRevE.98.042502}.

\bibitem{FranchiniBalzan2020}
S.~Franchini, R.~Balzan,
\textit{Energy of the interacting self-avoiding walk at the $\theta$ point},
Phys. Rev. E 102, 032143 (2020).
\url{https://doi.org/10.1103/PhysRevE.102.032143}.

\bibitem{Arthur1986}
W.~B.~Arthur, Yu.~M.~Ermoliev, Yu.~M.~Kaniovski,
\textit{Strong laws for a class of path-dependent stochastic processes with applications},
in \textit{Stochastic Optimization}, Lecture Notes in Control and
Information Sciences 81, 287--300 (1986).
\url{https://doi.org/10.1007/BFb0007105}.

\bibitem{Pemantle1991}
R.~Pemantle,
\textit{When are touchpoints limits for generalized P\'olya urns?},
Proc. Amer. Math. Soc. 113(1), 235--243 (1991).
\url{https://doi.org/10.1090/S0002-9939-1991-1055778-8}.

\bibitem{NakayamaMori2021}
K.~Nakayama, S.~Mori,
\textit{Universal function of the nonequilibrium phase transition of a nonlinear P\'olya urn},
Phys. Rev. E 104, 014109 (2021).
\url{https://doi.org/10.1103/PhysRevE.104.014109}.

\bibitem{HLS1987}
B.~M.~Hill, D.~Lane, W.~Sudderth,
\textit{Exchangeable urn processes},
Ann. Probab. 15(4), 1586--1592 (1987).
\url{https://doi.org/10.1214/aop/1176991995}.

\bibitem{Franchini_URNS}
S.~Franchini,
\textit{Large deviations for generalized Polya urns with arbitrary urn function},
Stochastic Process. Appl. 127(10), 3372--3411 (2017),
\url{https://doi.org/10.1016/j.spa.2017.02.010};
corrigendum: Stochastic Process. Appl. 189, 104745 (2025),
\url{https://doi.org/10.1016/j.spa.2025.104745}.

\bibitem{Jack2019}
R.~L.~Jack,
\textit{Large deviations in models of growing clusters with symmetry-breaking transitions},
Phys. Rev. E 100, 012140 (2019).
\url{https://doi.org/10.1103/PhysRevE.100.012140}.

\bibitem{BlumBrennan1980}
J.~R.~Blum, M.~Brennan,
\textit{On the strong law of large numbers for dependent random variables},
Israel J. Math. 37(3), 241--245 (1980).
\url{https://doi.org/10.1007/BF02760965}.

\bibitem{Franchini_K}
S.~Franchini,
\textit{Large deviations for Generalized Polya Urns with non-binary increments},
\href{https://arxiv.org/abs/2506.22234}{arXiv:2506.22234} (2026).
\url{https://doi.org/10.48550/arXiv.2506.22234}.

\bibitem{Kennedy2009}
T.~Kennedy,
\textit{Numerical computations for the Schramm--Loewner Evolution},
J. Stat. Phys. 137, 839--856 (2009).
\url{https://doi.org/10.1007/s10955-009-9866-2}.

\bibitem{Tran2015}
H.~Tran,
\textit{Convergence of an algorithm simulating Loewner curves},
Ann. Acad. Sci. Fenn. Math. 40, 601--616 (2015).
\url{https://doi.org/10.5186/aasfm.2015.4037}.

\bibitem{Bardella2024}
G.~Bardella, S.~Franchini, P.~Pani, S.~Ferraina,
\textit{Lattice physics approaches for neural networks},
iScience 27(12), 111390 (2024).
\url{https://doi.org/10.1016/j.isci.2024.111390}.

\bibitem{Friston2017}
K.~Friston, T.~FitzGerald, F.~Rigoli, P.~Schwartenbeck, G.~Pezzulo,
\textit{Active Inference: A Process Theory},
Neural Comput. 29(1), 1--49 (2017).
\url{https://doi.org/10.1162/NECO_a_00912}.

\bibitem{BandyopadhyayThacker2022}
A.~Bandyopadhyay, D.~Thacker,
\textit{A new approach to P\'olya urn schemes and its infinite color generalization},
Ann. Appl. Probab. 32(1), 46--79 (2022).
\url{https://doi.org/10.1214/21-AAP1671}.

\end{thebibliography}
\end{document}